\documentclass[10pt]{article}

\usepackage{latexsym}
\usepackage{amssymb}
\usepackage{epsfig}
\input amssym.def
\input amssym  
\usepackage{amssymb}
\newtheorem{theorem}{Theorem}[section]
\newtheorem{corollary}{Corollary}[section]
\newtheorem{lemma}{Lemma}[section]
\newtheorem{definition}[theorem]{Definition}

\newtheorem{conjecture}{Conjecture}[section]

\newtheorem{remark}{Remark}[section]
\newtheorem{proposition}{Proposition}[section]

\newcommand{\halmos}{\rule{1ex}{1.4ex}}
				
\newcommand{\bea}{\begin{eqnarray}}
\newcommand{\epf}{\hspace*{\fill}\mbox{$\halmos$}}
\newcommand{\eea}{\end{eqnarray}}
\newcommand{\nn}{\nonumber \\}
\newcommand{\be}{\begin {equation}}

\newcommand{\ee}{\end{equation}}

\title {{Vertex operator superalgebra structure for degenerate minimal
models: Neveu-Schwarz algebra}}
\author{Antun Milas \thanks{The author is partially supported by NSF grants.}}
\date{}
\begin{document}
\bibliographystyle{amsalpha}
\maketitle
\begin{abstract}
The $\mathbb{Z}/2\mathbb{Z}$--graded intertwining
operators are introduced. We study these
operators in the case of ``degenerate'' $N=1$ minimal models, with 
the central charge $c=\frac{3}{2}$. 
The corresponding fusion ring is isomorphic 
to the Grothendieck ring for the Lie superalgebra ${\mathfrak osp}(1|2)$.
We also discus multiplicity--$2$ fusion rules 
and logarithmic intertwiners.
\end{abstract}

\section{Introduction}

This paper is a continuation of
\cite{M1}. It is also closely related to \cite{HM}. 
For a more detailed introduction see \cite{M1}.

We introduce the notion of a $\mathbb{Z}/2\mathbb{Z}$--graded
({\em even} and {\em odd}) intertwining
operator and use it to study the {\em fusion ring} 
for the ``degenerate'' ($p=q$) minimal models, for the
vertex operator superalgebra $L(\frac{3}{2},0)$.
As in \cite{M1} there are two approaches: one which 
uses the lattice construction (extended with a
suitable fermionic Fock space) and the other which uses
the singular vectors and projection formulas.
For the future purposes we use the latter approach.

The degenerate minimal models are irreducible modules
for the $N=1$ superconformal vertex operator algebra $L(\frac{3}{2},0)$
(cf. \cite{KW}, \cite{Ad}) that are
isomorphic to $L(\frac{3}{2}, \frac{q^2}{2})$, $q \in \mathbb{N}$. 
We prove (see Theorem \ref{last} and Corollary \ref{last1})  that the corresponding
fusion ring is isomorphic to the Grothendieck ring
for ${\mathfrak osp}(1|2)$, i.e., we
formally have:
$$ L\left(\frac{3}{2},\frac{r^2}{2}\right) \times
L\left(\frac{3}{2},\frac{q^2}{2}\right)=$$
$$L\left(\frac{3}{2},\frac{(r+q)^2}{2}\right) +
L\left(\frac{3}{2},\frac{(r+q-1)^2}{2}\right)  \ldots
+L\left(\frac{3}{2},\frac{(r-q)^2}{2}\right),$$
for every $r,q \in {\bf N}$, $r \geq q$.

As in the Virasoro algebra case, these fusion coefficients
are $0$ or $1$. However in \cite{HM} we showed that
for some vertex operator algebras $L(c,0)$, fusion coefficients might be $2$.
In Proposition \ref{twons}
we constructe a non--trivial example when $c=\frac{15}{2}-3\sqrt{5}$.

At the very end, we constructe an example of a {\em logarithmic intertwining 
operator} (for the definition see \cite{M2})  for the $N=1$ vertex operator 
superalgebra  $L(\frac{27}{2},0)$ (cf. Proposition \ref{logns}).

\vskip 5mm 

\section{Lie superalgebra $ {\mathfrak osp}(1|2)$ and ${\mathcal
R}ep({\mathfrak osp}(2|1))$}
The Lie superalgebra ${\mathfrak osp}(1|2)$ is a graded extension of
the finite--dimensional Lie algebra ${\mathfrak sl}(2,{\bf C})$. 
It has three even generators $x,y$ and $h$, and two odd generators
$\varphi$ and $\chi$, that satisfy
$$[h,x]=2x, \ \ [h,y]=-2y, \ \ [x,y]=h,$$
$$[x,\chi]=\chi, \ \  [x,\varphi]=-\varphi, \ \ [y,\chi]=-\chi, \ \
[y,\varphi]=\varphi,$$
$$[h,\varphi]=-\varphi, \ \ [h,\chi]=\chi,$$
$$\{\chi,\varphi\}=2h, \{\chi,\chi\}=2x, \ \   \{\varphi,\varphi\}=2y .$$

Generators $\{x,y,h\}$ span a Lie algebra isomorphic to ${\mathfrak
sl}(2,\mathbb{C})$, and this fact
makes the representation theory of ${\mathfrak osp}(1|2)$ quite simple.
All irreducible $\mathfrak{osp}(1|2)$--modules can be constructed in the
following way. Fix a positive half integer $j$ ($2j \in {\bf N}$) and
a $4j+1$--dimensional vector
space $V(j)$ spanned by the vectors $\{v_{j},v_{j-1/2},...,v_{-j}\}$, with the
following actions: 
\begin{eqnarray}
&& x.v_{i}=\sqrt{[j-i][j+i+1]}v_{i+1}, \nn
&& y.v_{i}=\sqrt{[j+i][j-i+1]}v_{i-1}, \nn
&& h.v_{i}=2iv_{i}. 
\end{eqnarray}
If $2(i-j) \in {\bf Z}$ then we define
\begin{eqnarray}
&& \varphi.v_{i}=- \sqrt{j+i} v_{i-1/2}, \nn 
&& \chi.v_{i}=-\sqrt{j-i}v_{i+1/2},
\end{eqnarray}
otherwise
\begin{eqnarray}
&& \varphi.v_{i}=\sqrt{j-i+1/2} v_{i-1/2}, \nn 
&& \chi.v_{i}=-\sqrt{j+i+1/2}v_{i+1/2}.
\end{eqnarray}

In all these formulas $v_j=0$ if $j \notin \{j,j-\frac{1}{2},\ldots,-j \}$.
It is easy to see that each $V(j)$ is an irreducible $\mathfrak{osp}(1|2)$--module and
that every finite dimensional irreducible representation of $\mathfrak{osp}(1|2)$ is
isomorphic to some $V(j)$ for $j \in {\bf N}/2$.

It is a pleasant exercise to decompose tensor product $V(i) \otimes V(j)$. The following result is
well--known (see \cite{FM} for instance):
\begin{equation} \label{osp12}
V(i) \otimes V(j) \cong  \bigoplus_{k=|i-j|, k \in {\bf N}/2}^{i+j}
V(k).
\end{equation}

\section{$N=1$ Neveu-Schwarz superalgebra and its minimal models}

The $N=1$ Neveu-Schwarz superalgebra is given by
$$\mathfrak{n}\mathfrak{s}=\bigoplus_{n\in \mathbb{Z}}\mathbb{C}L_{n}\bigoplus
\bigoplus_{n\in \mathbb{Z}}\mathbb{C}G_{n+1/2}\bigoplus \mathbb{C}C,$$
together with the following
$N=1$ Neveu-Schwarz relations:
\begin{eqnarray*} \label{nscom}
{[L_{m}, L_{n}]}&=&(m-n)L_{m+n}+\frac{C}{12}
(m^{3}-m)\delta_{m+n, 0},\\
{[L_{m}, G_{n+ 1/2}]}&=&\left(\frac{m}{2}-\left(n+
\frac{1}{2}\right)\right)G_{m+n+ 1/2},\\
{[G_{m+1/2}, G_{n-1/2}]}&=& 
2L_{m+n}+\frac{C}{3}(m^{2}+m)\delta_{m+n, 0},\\
{[C, L_{m}]}&=&0,\\
{[C, G_{m+1/2}]}&=&0
\end{eqnarray*}
for $m, n\in \mathbb{Z}$. 
We have the standard triangular decomposition
$\mathfrak{ns}=\mathfrak{ns}_+ \oplus \mathfrak{ns}_0 \oplus \mathfrak{ns}_-$
(cf. \cite{HM}).
For every $(h,c) \in {\bf C}^2$, we denote  by
$M(c,h)$ Verma module for
$\mathfrak{ns}$ algebra.
For each $(p,q) \in {\bf N}^2$, $p = q \ {\rm mod}  \ 2$,
 let us introduce a family of complex 'curves'
$(h_{p,q}(t),c(t))$;
$$h_{p,q}(t)=\frac{1-p^2}{8}t^{-1}+\frac{1-pq}{4}+\frac{1-q^2}{8}t,$$
$$c(t)=\frac{15}{2}+3t^{-1}+3t.$$
Then from the determinant formula (see \cite{KWa} ) it follows
that $M(c,h)$ is reducible
if and only if there is a $t \in {\bf C}$ and $p,q \in {\bf N}$, $p= q  \ {\rm mod} \ 2$  such that
$c=c(t)$ and $h=h_{p,q}(t)$. In this case
$M(c,h)$ has a {\em singular vector} (i.e., a vector annihilated by
$\mathfrak{ns}_+$ ) of the weight $h+\frac{pq}{2}$.
Any such vector we denote by $v_{\frac{pq}{2}}$.

In this paper we are interested in the $t=-1$. Then
$c(-1)=\frac{3}{2}$ and
$h_{p,q}(-1)=\frac{(p-q)^2}{8}$.
$h_{p,q}(-1)=h_{1,p-q+1}(-1)$, so we consider only
the case $h_{1,q}:=h_{1,q}(-1)$, (here $q$ is odd and positive).
Hence, each Verma module $M(\frac{3}{2},h_{1,q})$ is
reducible.

The following result easily follows from \cite{D} (or \cite{Aa})
and \cite{KWa}:
\begin{proposition}
For every odd $q$, $M(\frac{3}{2},h_{1,q})$ has the
following embedding structure
\begin{equation}
\ldots \rightarrow M\left(\frac{3}{2},h_{1,q+4}\right) \rightarrow
M\left(\frac{3}{2},h_{1,q+2}\right) \rightarrow
M\left(\frac{3}{2},h_{1,q}\right)
\rightarrow 0.
\end{equation}
Moreover, we have the following exact sequence:
\begin{equation}
0 \rightarrow M\left(\frac{3}{2},h_{1,q+2}\right) \rightarrow
M\left(\frac{3}{2},h_{1,q}\right) \rightarrow
L\left(\frac{3}{2},h_{1,q}\right) \rightarrow 0,
\end{equation}
where $L(\frac{3}{2},h_{1,q})$ is the corresponding
irreducible quotient.
\end{proposition}
\epf

Benoit and Saint-Aubin (cf. \cite{BSA}) found an explicit expression
for the singular vector $v_{1,q} \in M
\left(\frac{3}{2},h_{1,q}\right)$ that generates 
the maximal submodule:
\begin{equation}
v_{1,q}=\sum_{N;k_1,...,k_N} \sum_{\sigma \in S_N}  (-1)^{\frac{q-N}{2}}
c(k_{\sigma(1)},...,k_{\sigma(k)})
G_{-k_1/2}\ldots G_{-k_N/2}v,
\end{equation}
where $S_N$ is a symmetric group on $N$ letters
and the first summation is over all the partitions of $q$ into the odd
integers $k_1,..,k_N$ and
$$c(k_{\sigma(1)},...,k_{\sigma(k)})=\prod_{i=1}^N {k_i-1 \choose
(k_i-1)/2 } \prod_{j=1}^{(N-1)/2} \frac{4}{\sigma_{2j} \rho_{2j}},$$
where $\sigma_j=\sum_{l=1}^j k_l$ and $\rho_j=\sum_{l=j}^N k_l$.

In the special case: $q=1$, $h_{1,1}=0$, 
$M(\frac{3}{2},0)$ has a singular vector
$G(-1/2)v$ which generate the maximal submodule.
By quotienting we obtain a {\em vacuum} module
$L(\frac{3}{2},0)=M(\frac{3}{2},0)/\langle G_{-1/2}v_{3/2,0} \rangle$.

\section{$N=1$ superconformal vertex operator superalgebra and intertwining 
operators}

We use the definition of $N=1$ superconformal vertex operator superalgebra
(with and without odd variables) as 
in \cite{B} and \cite{HM} (see also \cite{KW} and \cite{Kac}).

Let $\varphi$ be a Grassman (odd) variable such that $\varphi^2=0$.
Every $N=1$ superconformal vertex operator superalgebra 
$(V, Y, \mathbf{1}, \tau)$ 
can be equipped with a structure of  $N=1$ superconformal vertex operator algebra
with odd variables via
\begin{eqnarray*}
Y( \ ,(x,\varphi)) : V\otimes V&\to &V((x))[\varphi], \\
u\otimes v&\mapsto& Y(u, (x, \varphi)) v,
\end{eqnarray*}
where
$$Y(u, (x, \varphi))v=Y(u, x)v+\varphi Y(G(-1/2)u, x)v$$
for $u, v\in V$.

The same formula can be used in the case of modules
for the superconformal vertex operator superalgebra $(V, Y,
\mathbf{1}, \tau)$ (see \cite{HM}).

It is known (\cite{KW}) that 
$V(c,0):=M(c,0)/ \langle G_{-1/2}v_{c,0} \rangle$ 
\footnote{We write $L(c,0)$ if $V(c,0)$ is
irreducible.} is
a $N=1$ superconformal vertex operator superalgebra.
Also, every lowest weight $\mathfrak{ns}$--module with the central charge
$c$, is a $V(c,0)$--module.
If $c=\frac{3}{2}$ then $V(\frac{3}{2},0)=L(\frac{3}{2},0)$. Hence
\begin{proposition}
Every irreducible
$L(\frac{3}{2},0)$--module
is isomorphic to 
$L(\frac{3}{2},h)$, for some $h \in {\bf C}$.
\end{proposition}
{\em Proof:} It is known (cf. \cite{KW}) that there is one--to--one
correspondence between equivalence clasess of irreducible $A(L(\frac{3}{2},0))$--modules (here
$A(L(\frac{3}{2},0)) \cong \mathbb{C}[y]$ is Zhu's associative algebra) and
irreducible $L(\frac{3}{2},0)$--modules. If $W$ is an irreducible
$A(L(\frac{3}{2},0))$--module, then there is an irreducible
$L(\frac{3}{2},0)$--module $\Omega(W)$  (that is $\mathbb{N}$--gradable) 
such that $\Omega(W)(0) \cong W$. Because $\Omega(W)$ is
a $\goth{ns}$--irreducible module it is also
$\mathbb{N}$--gradable.
Therefore $\Omega(W) \cong L(c,h)$ for some $h \in \mathbb{C}$.
\epf

Among all irreducible $L(\frac{3}{2},0)$--modules
we distinguish modules isomorphic to $L(\frac{3}{2},h_{1,q})$, 
$ q \in 2{\bf N}-1$. These are so called degenerate minimal models.

\subsection{Intertwining operators and its matrix coefficients}

The notation of an intertwining operators for $N=1$ superconformal vertex
operator algebras is introduced in \cite{KW} and \cite{HM}.

Let $W_{1}$, $W_{2}$ and $W_{3}$ be a triple of $V$--modules
and $\mathcal{Y}$ an intertwining operator of type ${W_3
\choose W_{1} \ W_{2}}$. 
Then we consider the corresponding intertwining operator   
with odd variable (cf. \cite{HM}):
\begin{eqnarray*}
\mathcal{Y}( \ , (x,\varphi)) \ : W_{1}\otimes W_{2}&\to &W_{3}\{x\}[\varphi]\\
w_{(1)}\otimes w_{(2)}&\mapsto& \mathcal{Y}(w_{(1)}, (x, \varphi)) w_{(2)},
\end{eqnarray*}
such that 
$$\mathcal{Y}(w_{(1)}, (x, \varphi)) w_{(2)}
=\mathcal{Y}(w_{(1)}, x) w_{(2)}+
\varphi\mathcal{Y}(G(-1/2)w_{(1)}, x) w_{(2)}.$$

Let $w_1$ be a lowest weight vector for the Neveu-Schwarz algebra of
the weight $h$.
From the Jacobi identity we derive the following formulas:
\bea \label{evencom}
&& [L(-n),{\mathcal Y}(w_1,x_2)]=(x_2^{-n+1}\frac{\partial}{\partial
x_2}+(1-n)h){\mathcal Y}(w_1,x_2), \nn
&&  [G(-n-1/2),{\mathcal Y}(w_1,x_2)]=x_2^{-n}{\mathcal Y}(G(-1/2) w_1,x_2), \nn
&& [L(-n),{\mathcal Y}(G(-1/2)w_1,x_2)]=(x_2^{-n+1}\frac{\partial}{\partial
x_2}+(1-n)(h+\frac{1}{2}){\mathcal Y}(G(-1/2)w_1,x_2), \nn
&&  [G(-n-1/2),{\mathcal Y}(G(-1/2)w_1,x_2)]=
(x_2^{-n}\frac{\partial}{\partial x_2}-2nhx_2^{-n-1}){\mathcal
Y}(w_1,x_2).
\eea
In the odd formulation we obtain
\begin{eqnarray} \label{oddcom}
&& [L(-n),\mathcal{ Y}(w_{1},(x_2,\varphi))] \nn
&&= (x_2^{-n+1}\partial_{x_2}+(1-n)x_2^{-n}(h
+1/2\varphi \partial_{\varphi}))\mathcal{ Y}(w_{1},(x_2,\varphi)) \nn
&& [G(-n-1/2),\mathcal{Y}(w_{1},(x_2,\varphi))] \nn
&&= (x_2^{-n}(\partial_{\varphi}-\varphi \partial_{x_2})-2n
x_2^{-n-1}(h \varphi)\mathcal{ Y}(w_{1},(x_2,\varphi)),
\end{eqnarray}
where $\partial_{\varphi}$ is the odd (Grassmann) derivative.

\subsection{Even and odd intertwining operators}

In \cite{HM} we proved that every intertwining operator
$$\mathcal{Y} \in I {L(c,h_3) \choose L(c,h_1) \ L(c,h_2)}$$
is uniquely determined by the operators ${\mathcal Y}(w_1,x)$
and ${\mathcal Y}(G(-1/2)w_1,x)$, where $w_1$ is the highest weight
vector of $L(c,h_1)$. This fact will be used later in connection 
with the following definition. 

\begin{definition} \label{oddeven}
{\em Let $| \ |$ denote the parity ($0$ or $1$).
An intertwining operator ${\mathcal Y} \in I {W_3 \choose W_1 \ W_2}$ is:
\begin{itemize}
\item {\em even}, if 
$${\rm Coeff}_{x^s}|{\mathcal Y}(w_1,x)w_2|=|w_1|+|w_2|,$$
\item {\em odd}, if
$${\rm Coeff}_{x^s}|{\mathcal Y}(w_1,x)w_2|=|w_1|+|w_2|+1,$$
\end{itemize}
for every $s \in \mathbb{C}$ and 
every $\mathbb{Z}/2\mathbb{Z}$--homogeneous vectors $w_1$ and $w_2$.}
\end{definition}

The space of even (odd) intertwining operators of the type
${\mathcal Y} \in I {W_3 \choose W_1 \ W_2}$ we denote
by $I \ {W_3 \choose W_1 \ W_2}_{\rm even}$ ($I \ {W_3 \choose W_1 \ W_2}_{\rm
odd}$).
In general one does not expect a decomposition of 
$I {W_3 \choose W_1 \ W_2}$ into the even and the odd subspace.

\subsection{Frenkel-Zhu's theorem for vertex operator superalgebras}
According to \cite{KW} (after \cite{Z}),
to every vertex operator superalgebra we can associate 
Zhu's associative algebra $A(V)$. 
If $V=L(c,0)$, $A(L(c,0)) \cong {\bf C}[y]$.
where $y=[(L(-2)-L(-1)){\bf 1}]=[L(-2){\bf 1}]$ 
(because of the calculations that follow it is convenient to use
$y=[(L(-2)-L(-1)){\bf 1}]$).
Also to every $V$--module 
$W$  we associate a $A(V)$--bimodule $A(W)$ (cf. \cite{KW}).
In a special case $W=M(c,h)$, we have
$$A(M_{\mathfrak{ns}}(c,h))=M_{\mathfrak{ns}}(c,h)/O(M_{\mathfrak{ns}}(c,h)),$$
where 
\begin{eqnarray} \label{bimodules}
&& O(M_{\mathfrak{ns}}(c,h))=  \{L(-n-3)-2L(-n-2)+L(-1)v, \nn 
&& G(-n-1/2)-G(-n-3/2)v, n \geq 0, v \in M(c,h) \}.
\end{eqnarray}
It is not hard to see that, as $\mathbb{C}[y]$--bimodule, 
$$A(M(c,h)) \cong {\bf C}[x,y] \oplus {\bf C}[x,y]v,$$
where $v=[G(-1/2)v_{h}]$ and
$$y=[L(-2)-L(-1)], \ \ x=[L(-2)-2L(-1)+L(0)].$$

Let $W_{1}$, $W_{2}$ and $W_{3}$ be three 
$\mathbb{N}/2$--gradable irreducible $V$--modules
such that ${\rm Spec}L(0)|_{W_i} \in h_i + \mathbb{N}$, $i=1,2,3$
and $\mathcal{Y} \in I \ {W_3 \choose W_1 \ W_2}$. We define
$o(w_1):={\rm Coeff}_{x^{h_3-h_1-h_2}}\mathcal{Y}(w_1,x)$.
Because the fusion rules formula in \cite{FZ}
needs some modifications (cf. \cite{L1}) the
same modification is necessary for the main Theorem in \cite{KW}
(this can be done with a minor super--modifications along the lines of
\cite{L1}). Nevertheless (cf. \cite{KW}):
\begin{theorem} \label{superfz}
The mapping 
$$\pi : I {W_3 \choose W_{1} \ W_{2}} \rightarrow
{\rm Hom}_{A(V)}(A(W_{1})\otimes _{A(V)}W_{2}(0), W_{3}(0)),$$
such that
\begin{equation}
\pi(\mathcal{Y})(w_1 \otimes w_2)=o(w_1)w_2,
\end{equation}
is injective.
\end{theorem}

\section{Some Lie superalgebra homology}

In this section we recall
some basic definition from the homology theory
of infinite dimensional Lie superalgebras
which is in the scope of the monograph \cite{F} (in the cohomology
setting though).

Let $\mathcal L$ be an any (possibly infinite dimensional) $\mathbb
{Z}/ 2\mathbb{Z}$--graded Lie
superalgebra with the $\mathbb{Z}/2\mathbb{Z}$--decomposition:
${\mathcal L}={\mathcal L}_0 \oplus {\mathcal L}_1$. 
and let $M=M_0 \oplus M_1$ be any ${\bf Z}_2$--graded
${\mathcal L}$--module, such that the gradings are compatible.
Then, we form a chain complex $(C,d,M)$ (for details see \cite{F}),
$$0 \stackrel{d_0}{\leftarrow} C_0({\mathcal L},M) \stackrel{d_1}{\leftarrow}
C_1({\mathcal L},M) 
\stackrel{d}{\leftarrow }\ldots ,$$
where 
$$C_q({\mathcal L},M)=\bigoplus_{q_0+q_1=q} M \otimes \Lambda^{q_0}{\mathcal L}_0 \otimes
S^{q_1}{\mathcal L}_1,$$
$$C_q^p({\mathcal L},M)=\bigoplus_{\stackrel{q_0+q_1=q }{q_1+r = p \ {\rm
mod} 2 } } 
M_r \otimes \Lambda^{q_0}{\mathcal L}_0 \otimes S^{q_1}{\mathcal L}_1,$$
for $p=0,1$.
The mappings $d$ are super--differentials.
For $q \in {\bf N}$ and $p=0,1$, we define $q$--th homology
with coefficients in $M$ as:
\begin{equation} \label{z2}
H_q^p({\mathcal L}, M)=
{\rm Ker}(d_q (C_q^p({\mathcal L},M)))_p/(d_{r+1}(C_{q+1}^p({\mathcal
L},M)))_p.
\end{equation}
In a special case $q=0$, we have
$$H_0^0({\mathcal L}, M)=M_0/({\mathcal L}_0 M_0 + {\mathcal L}_1 M_1),$$
and
$$H_0^1({\mathcal L}, M)=M_1/({\mathcal L}_1 M_0 + {\mathcal L}_0 M_1).$$

We want to calculate $H_q({\mathcal L}_{s},L(\frac{3}{2},h_{1,q}))$.
for the Lie superalgebra
$${\mathcal L}_s = \bigoplus_{n \geq 0} {\mathcal L}_s(n),$$
where $ {\mathcal L}_s(n)$ is spanned by the vectors
$L(-n-3)-2L(-n-2)+L(-n-1)$ and $G(-n-1/2)-G(-n-3/2)$
, $n \in {\bf N}$.
From (\ref{bimodules}) we see (cf. \cite{HM}) 
that $H_0({\mathcal L}_{s},M(c,h))$ is a 
$\mathbb{C}[y]$--bimodule such that:
\begin{equation}
H_0({\mathcal L}_{s},M(c,h)) \cong A(M(c,h)) \cong {\bf C}[x,y] \oplus
{\bf C}[x,y]v.
\end{equation} 

\begin{remark}
{\em It is more involved to calculate $H_0(({\mathcal L}_{s},L(c,h))$, so we
consider only the special case $c=\frac{3}{2}$, $h=h_{1,q}$,
$q$ odd. As in the Virasoro case (see \cite{M1}) it is easy to show that
the space $H_p({\mathcal L}_{s},L(\frac{3}{2},h_{1,q}))$ is infinite
dimensional for very $p,q,s \in {\bf N}$, and finitely generated as
a $A(L(3/2,0))$--module.
Moreover,
$${\rm Ext}^1(L(\frac{3}{2},h_{1,q}),L(\frac{3}{2},h_{1,r}))$$
is one-dimensional if $r=q+2$ and $0$ otherwise.
We will not need these results.}
\end{remark}

In the minimal models case we expect 
a substantially different result (cf \cite{FF1}). 
\begin{conjecture} \label{conj}
Let $c_{p,q}=\frac{3}{2}\left( 1-2\frac{(p-q)^2}{pq}\right)$ and
$h_{p,q}^{m,n}=\frac{(np-mq)^{2}-(p-q)^{2}}{8pq}$. Then
$${\rm dim} \ H_{q}({\mathcal L}_s, L(c_{p,q},h_{p,q}^{m,n}))< \infty,$$
for every $q \in {\bf N}$.
\end{conjecture}
There is strong evidence that Conjecture (\ref{conj})
holds based on \cite{Ad} and an example $c=-\frac{11}{14}$ treated in
Appendix of \cite{HM}.

The main difference between the
minimal models and the degenerate models is the fact that the maximal 
submodule for a minimal model is generated by two singular 
vectors, compared to 
$M(\frac{3}{2},h_{1,q})$ where the maximal submodule
is generated by a single singular vector.

\section{Benoit-Saint-Aubin's formula projection formulas} 
\subsection{Odd variable formulation}
We have seen before how to derive the commutation relation between
generators of $\mathfrak{ns}$ superalgebra and
${\mathcal Y}(w_1,x)$ where $w_1$ is a lowest weight vector for $\goth{ns}$.
We fix ${\mathcal Y} \in I {L(\frac{3}{2},h) \choose
L(\frac{3}{2},h_{1,r}) \ L(\frac{3}{2},h_{1,q})}$ 
and consider the following matrix coefficient, 
\begin{equation}
\langle w'_{3}, \mathcal{Y}(w_1,x,\varphi)P_{{\rm sing}} w_2 \rangle,
\end{equation}
where $P_{sing}w_2=v_{1,q}$ (${\rm deg}(P_{sing})=q/2$) 
and $w_i$, $i=1,2,3$ are the lowest weight vectors. 

Since all modules are irreducible, by using
a result from \cite{HM} (Proposition 2.2), we get
$$\langle w'_{3}, \mathcal{Y}(w_1,x,\varphi)w_2 \rangle=c_1
x^{h-h_{1,q}-h_{1,r}} + c_2  \varphi x^{h-h_{1,q}-h_{1,r}-1/2},$$
where $c_1$ and $c_2$ are constants with the property 
\bea
&& c_1=c_2=0 \ {\rm implies} \ \ {\mathcal Y}=0.
\eea
From the formulas (\ref{oddcom})
$$\langle w'_{3}, \mathcal{Y}(w_1,x,\varphi)P_{{\rm sing}} w_2 \rangle
=P(\partial_{x_2},\varphi) \langle w'_{3},
\mathcal{Y}(w_1,x,\varphi)w_2 \rangle,$$
where $P(\partial_{x_2},\varphi)$ is a certain
super-differential operator such that
$${\rm deg}(P_{{\rm sing}})={\rm deg}P(\partial_{x_2},\varphi)=q/2.$$
Therefore
$$P(\partial_{x_2},\varphi)c_1x^{h-h_{1,q}-h_{1,r}}=\varphi
C_1(h_{1,q},h_{1,r},h)x^{h-h_{1,q}-h_{1,r}-q/2},$$
and
$$P(\partial_{x_2},\varphi)\varphi c_2x^{h-h_{1,q}-h_{1,r}-q/2}=
C_2(h_{1,q},h_{1,r},h)x^{h-h_{1,q}-h_{1,r}-q/2}.$$

Constants $C_1(h_{1,q},h_{1,r},h)$ and 
$C_2(h_{1,q},h_{1,r}, h)$ (in slightly different form, but
in more general setting) were
derived in \cite{BSA}. 
Considering these
coefficients was motivated by deriving formulas for singular vectors from
already known singular vectors. 
By slightly modifying result from \cite{BSA} we obtain 
\begin{proposition} \label{main}
Suppose that ${\mathcal Y} \in I {L(\frac{3}{2},h) \choose
L(\frac{3}{2},h_{1,r}) \ L(\frac{3}{2},h_{1,q})}$ and
$P(\partial_x,\varphi)$ are as the above 
Then, up to a multiplicative constant,
$$C_1(h_{1,q},h_{1,r},h)=\prod_{-j \leq k \leq j} (h-h_{1,q+4k})$$
and
$$C_2(h_{1,q},h_{1,r},h)=\prod_{-j+1/2 \leq k \leq j-1/2}
(h+\frac{1}{2}-h_{1,q+4k}),$$
for $j=(r-1)/4$, $j >0$ (when $j=0$, $C_2(h_{1,1},h_{1,r},h)=1$).
\end{proposition}
{\em Proof:} 
The superdifferential operator $P(\partial_x,\varphi)$
is obtained by replacing generators $L(-m)$ and $G(-n-1/2)$
by the superdifferential operators
\begin{equation} \label{lm}
L(-m) \mapsto -(x_2^{-m+1}\partial_{x_2}+(1-m)x_2^{-m}(h_1
+1/2\varphi \partial_{\varphi}))
\end{equation}
and
\begin{equation}
G(-n-1/2) \mapsto (x_2^{-n}(\partial_{\varphi}-\varphi \partial_{x_2})-2n
x_2^{-n-1}(h_1 \varphi)),
\end{equation}
acting on $\langle w'_3, \mathcal{Y}(w_1,x,\varphi)w_2 \rangle$.
This action was calculated in \cite{BSA}. 
Their results (Formula 3.10 in
\cite{BSA}) implies the statement \footnote{In \cite{BSA}
a different sign was used in the equation (\ref{lm}). Still, we obtain
the same result if we consider an isomorphic algebra with the generators
$\tilde{L}(n):=-L(n)$. The same generators were used in \cite{FF2}.}.
\epf

\subsection{BSA formula without odd variables}
Since Frenkel-Zhu's formula does not involve odd variables
we need a version of Proposition \ref{main} without odd variables (which is of
course equivalent).
Again 
${\mathcal Y} \in  I \ {L(3/2,h) \choose L(3/2,h_{1,r}) \
L(3/2,h_{1,q})}$ is the same as the above. 
Then
$$\langle w'_{3}, \mathcal{Y}(w_1,x)P_{sing} w_2 \rangle=P_2(\partial_{x})
\langle w'_{3}, \mathcal{Y}(G(-1/2)w_1,x)w_2 \rangle,$$
and
$$\langle w'_{3}, \mathcal{Y}(G(-1/2)w_1,x)P_{sing} w_2 \rangle=P_1(\partial_{x})
\langle w'_{3}, \mathcal{Y}(w_1,x)w_2 \rangle,$$
where $P_1$ and $P_2$ are certain differential operators.
If 
$$P_2(\partial_{x})c_2x^{h-h_{1,q}-h_{1,r}-1/2}=c_2 K_2(h_{1,q}, h_{1,r},
h)x^{h-h_{1,q}-h_{1,r}-q/2 },$$
and 
$$P_1(\partial_{x})c_1x^{h-h_{1,q}-h_{1,r}}=c_1 K_1(h_{1,q}, h_{1,r}, h )
x^{h-h_{1,q}-h_{1,r}-q/2},$$
then, by comparing corresponding coefficients, we obtain 
\begin{eqnarray} \label{matching}
&& K_1(h_{1,q}, h_{1,r}, h)=C_1(h_{1,q}, h_{1,r}, h), \nn
&& K_2(h_{1,q}, h_{1,r}, h)=C_2(h_{1,q}, h_{1,r}, h).
\end{eqnarray}

Let us mention that the projection formulas from Proposition \ref{main}
have a simple explanation in the term of {\em super density modules}
for the Neveu-Schwarz superalgebra.

\section{Fusion ring for the degenerate minimal models}

In order to obtain an upper bound for the fusion
coefficients (cf. Theorem \ref{superfz})
we first compute 
$$A(L(\frac{3}{2},h_{1,q})) \otimes_{A(L(3/2,0)}
L(\frac{3}{2},h_{1,r})(0).$$

${\mathbb Z}_2$--grading of the $0$--th homology group (\ref{z2}) will enable
us (see Theorem (\ref{last}) to study odd and even intertwining
operators (see Definition \ref{oddeven}).
For that purpose we introduce the following splitting:
 
\begin{eqnarray} \label{zeroth}
&& A^0(L(\frac{3}{2},h_{1,q})):=H^0_0 ({\mathcal L}_s, L(\frac{3}{2},h_{1,q}))
\cong \frac{{\bf C}[x,y]}{I_1} \nn
&& A^1(L(\frac{3}{2},h_{1,q})):=H^1_0 ({\mathcal L}_s, L(\frac{3}{2},h_{1,q}))
\cong \frac{{\bf C}[x,y]v}{I_2},
\end{eqnarray}
where $I_1$ and $I_2$ are cyclic submodules (the maximal submodule
for $M(3/2,h_{1,q})$ is cyclic !).
It seems hard to obtain explicitly these polynomials.
First we obtain some useful formulas
Inside $A(M(c,h))$ (cf. \cite{W}):
\bea \label{ok1}
&& [L(-n)v]=[((n-1)(L(-2)-L(-1))+L(-1))v]=\nn
&& [(n(L(-2)-L(-1))-(L(-2)-2L(-1)+L(0))+L(0))v]=\nn
&& (ny-x+{\rm wt}(v))[v].
\eea
for every $n \in \bf{N}$ and every homogeneous $v \in M(c,h)$.
Therefore in
$$A(M(\frac{3}{2},h_{1,q}))
\otimes_{A(L(\frac{3}{2},0))} L(\frac{3}{2},h_{1,r})(0)$$
we have
\begin{eqnarray} \label{ok32}
&& [L(-n)v]=(nh_{1,q}-x+L(0))[v]. \nn
&& [G(-n-1/2)v]=[G(-1/2)v].
\end{eqnarray}
Also, we have:
\bea \label{ok2}
&& [G(-n-\frac{1}{2})G(-m-\frac{1}{2})v]=[G(-1/2)G(-m-1/2)v]=\nn
&& [(2L(-m-1)-G(-m-1/2)G(-1/2))v]=[(2L(-m-1)-L(-1))v]=\nn
&& ((2m+1)y-x+{\rm wt}(v))[v].
\eea
By using (\ref{ok1}) and (\ref{ok2}) we obtain
\bea \label{ok3}
&& [G(-m_1-1/2) \ldots G(-m_{2r}-1/2)L(-n_1)...L(-n_s)v_{1,q}]=\nn
&& \prod_{i=1}^r ((2m_{2i}+1)h_{1,r}-x+\sum_{p=2i+1}^{2r}(m_p+1/2)+h_{1,q})
\cdot \nn
&& \prod_{j=1}^s (n_j h_{1,r}-x+\sum_{p=j+1}^s n_p+h_{1,q})[v].
\eea
inside
$$A(M(\frac{3}{2},h_{1,q}))
\otimes_{A(L(\frac{3}{2},0))} L(\frac{3}{2},h_{1,r})(0).$$
It is easy to obtain a similar formula for the vector
$$[G(-m_1-1/2) \ldots G(-m_{2r+1}-1/2)L(-n_1) \cdots L(-n_s)v_{1,q}].$$

\begin{lemma} \label{lelast}
Let $[P_{\rm sing}v_{1,q}]=Q_1(x)[G(-1/2)v_{1,q}]$ and
$[G(-1/2)P_{\rm sing}v_{1,q}]=Q_2(x)[v_{1,q}]$ be projections inside
$$A(M(\frac{3}{2},h_{1,q}))
\otimes_{A(L(\frac{3}{2},0)} L(\frac{3}{2},h_{1,r})(0).$$
Then 
\bea
&& Q_1(h)=K_2(h_{1,q},h_{1,r},h), \nn
&& Q_2(h)=K_1(h_{1,q},h_{1,r},h),
\eea
for every $h \in \mathbb{C}$.
\end{lemma}
{\em Proof:} 
We use the notation from the section 6.2, where
$$\mathcal{Y} \in I \ {L(3/2,h) \choose L(3/2,h_{1,r}) \ L(3/2,h_{1,q})}.$$ 
By using (\ref{evencom}), we obtain
\bea \label{ok4}
&& \langle w'_3,
{\mathcal Y}(w_1,x)G(-m_1-1/2)...G(-m_{2r}-1/2)L(-n_1)...L(-n_s)w_2
\rangle= \nn
&& \prod_{i=1}^r -(x^{-m_{2i-1}-m_{2i}}\frac{\partial}{\partial x}
-2m_{2i}h_{1,r}x^{-m_{2i-1}-m_{2i}-1}) \cdot \nn
&& \prod_{j=1}^s
-(x^{-n_j+1}\frac{\partial}{\partial x}+(1-n_j)h_{1,r}x^{-n_j}) 
\langle w'_3, {\mathcal Y}(w_1,x)w_2 \rangle=\nn
&& c_1 \prod_{i=1}^r ((2m_{2i}+1)h_{1,r}-h+h_{1,q}+\sum_{p=2i+1}^{2r}
(m_{p}+1/2)) \cdot \nn
&& \prod_{j=1}^s (n_j h_{1,r}-h+\sum_{p=j+1}^s
n_p+h_{1,q})x^{h-h_{1,q}-
h_{1,r}-r-\sum m_i -\sum_j n_j},
\eea 
for the constant $c_1$ (see Section 6.1 and 6.2) that depends only on $\mathcal{Y}$.
There is a similar expression for
\begin{equation} \label{ok5}
\langle w'_3,
{\mathcal Y}(w_1,x)G(-m_1-1/2)...G(-m_{2r+1}-1/2)L(-n_1)...L(-n_s)w_2
\rangle.
\end{equation}
If we compare (\ref{ok3}) with (\ref{ok4}) (and corresponding
formulas for (\ref{ok5})) it follows that
$Q_1(h)$ is, up to a non--zero multiplicative constant, 
equal to $K_2(h_{1,r},h_{1,q},h)$ (singular vector is odd!)
and $Q_2(h)$ is, up to a multiplicative constant, equal to
$K_1(h_{1,r},h_{1,q},h)$.
\epf

Thus, Proposition \ref{main}
and Theorem \ref{last} gives us
\begin{theorem} \label{last}
\begin{itemize}
\item[(a)]
As a $A(L(3/2,0))$--module \\
\bea \label{amdec}
&& A(L(3/2,h_{1,q})) \otimes_{A(L(3/2,0))} L(3/2,h_{1,r})(0)
\cong \nn
&& \frac{{\bf C}[x]}{<\prod_{-j \leq k \leq j}(x-h_{1,q+4k})>} \oplus 
\frac{{\bf C}[x]}{<\prod_{-j+1/2 \leq k \leq j+1/2} (h+1/2-h_{1,q+4k})>}.
\eea
\item[(b)]
The space 
$$I \ {L(3/2,h) \choose L(3/2,h_{1,q}) \ M(3/2,h_{1,r})},$$
is non--trivial if and only if $h=h_{1,s}$ for
some $s \in \{q+r-1,q+r-3,\ldots,q-r+1\}$.
\item[(c)]
The space 
$$I \ {L(3/2,h) \choose L(3/2,h_{1,q}) \ L(3/2,h_{1,r})},$$
is one--dimensional if and only if $h=h_{1,s}$, $s \in
\{q+r-1,q+r-3,\ldots,|q-r|+1\}.$
\end{itemize}
\end{theorem}
{\em Proof (a):}
From  Lemma \ref{lelast} it follows that
\bea
&& A(L(3/2,h_{1,r})) \otimes_{A(L(3/2,0))} L(3/2,h_{1,q}) \cong 
\frac{\mathbb{C}[x]}{\langle  Q_1(x) \rangle} \oplus \frac{\mathbb{C}[x]}{\langle  Q_2(x) \rangle}.
\eea
Now we apply (\ref{matching}) and Proposition \ref{main}.\\
{\em Proof (b):}
As in \cite{M1}, by examining carefully the main construction of 
intertwining operators in \cite{L1} with a minor super--modifications,
for every $A(L(3/2,0))$--morphism from 
$A(L(3/2,h_{1,q})) \otimes_{A(L(3/2,0))} L(3/2,h_{1,r})$ to $L(3/2,h)(0)$
we can construct a non--trivial intertwining operator of the form
$I \ {L(3/2,h) \choose L(3/2,h_{1,q}) \ M(3/2,h_{1,r})}.$ \\
{\em Proof (c):} In order to project 
$\mathcal{Y} \in I \ {L(3/2,h) \choose L(3/2,h_{1,q}) \
M(3/2,h_{1,r})}$ to a non--trivial intertwining operator of the type
$ {L(3/2,h) \choose L(3/2,h_{1,q}) \ L(3/2,h_{1,r})}$
(as in \cite{M1}) $h=h_{1,s}$ for 
$s \in \{q+r-1,q+r-3,\ldots,q-r+1\} \cap
\{r+q-1,r+q-3,\ldots,r-q+1\}$, i.e.,
$s \in \{q+r-1,q+r-3,\ldots,|q-r|+1\}$.
\epf

\begin{theorem}
Suppose that $q \geq r$ \footnote{${W_3 \choose W_1 \ W_2} \cong 
{W_3 \choose W_2 \ W_1}$.}
\begin{equation}
{\rm dim} \ I \ {L(3/2,h_{1,s}) \choose L(3/2,h_{1,q}) \
L(3/2,h_{1,r})}_{\rm even}=1,
\end{equation}
if and only if
$$s \in \{ q+r-1,q+r-5,...,q-r+1 \}$$
\begin{equation}
{\rm dim} \ I \ {L(3/2,h_{1,s}) \choose L(3/2,h_{1,q}) \
L(3/2,h_{1,r})}_{\rm odd}=1
\end{equation}
if and only if
$$s \in \{ q+r-3,q+r-7,...,q-r+3 \}.$$
\end{theorem}
{\em Proof:}
By using (\ref{amdec}) we obtain 
the following decomposition:
\begin{eqnarray} \label{decomposition}
&& A^0(L(\frac{3}{2},h_{1,q})) \otimes_{A(L(3/2,0))} L(\frac{3}{2},h_{1,r})(0)
\cong \nn
&& {\mathbb C} v_{q+r-1} \oplus {\mathbb C} v_{q+r-5} \ldots \oplus {\mathbb C}v_{q-r+1}
\nn
&& A^1 L(\frac{3}{2},h_{1,q})) \otimes_{A(L(3/2,0))} L(\frac{3}{2},h_{1,r})(0)
\cong \nn
&& {\mathbb C} v_{q+r-3} \oplus {\mathbb C} v_{q+r-7} \ldots \oplus {\mathbb C}v_{q-r+3},
\end{eqnarray}
where ${\mathbb C}v_{i}$ is a ${\mathbb C}[y]$--module such that
$y.v_i=\frac{(i-1)^2}{8}v_i$. \\
{\em Claim:} Let 
$$\psi \in {\rm Hom}_{A(L(c,0))}(A^0(L(\frac{3}{2},h_{1,q}))
\otimes_{A(L(3/2,0))}L(\frac{3}{2},h_{1,r})(0),
L(\frac{3}{2},h_{1,s})(0)),$$
then the corresponding intertwining operator is even. 
Similarly if we start from
$$\psi \in {\rm Hom}_{A(L(c,0))}(A^1(L(\frac{3}{2},h_{1,q}))
\otimes_{A(L(3/2,0))}L(\frac{3}{2},h_{1,r})(0),
L(\frac{3}{2},h_{1,s})(0)),$$
the corresponding intertwining operator is odd.
\\
{\em Proof (Claim):} Let us elaborate the proof when $\psi$ is ``even''.
From the construction in \cite{FZ} and \cite{L2}
$\mathcal{Y}$ is obtained by lifting $\psi$ to a mapping from
$L(3/2,h_{1,q}) \otimes L(3/2,h_{1,r})(0)$ to
$L(3/2,h_{1,s})(0)$, such that
$$L(3/2,h_{1,q})_{\rm odd} \otimes L(3/2,h_{1,r})(0) \mapsto 0.$$
To extend this map to a mapping from
$L(3/2,h_{1,q}) \otimes M(3/2,h_{1,r})$ to $L(3/2,h_{1,s})$
one uses generators and PBW so the sign is preserved. The last step
(projection to $L(3/2,h_{1,q}) \otimes L(3/2,h_{1,r})$) is possible
(because of the condition $q \geq r$) so the proof follows (when
$\psi$ is odd a similar argument works).
\epf

Let us summarize everything.
\begin{corollary} \label{last1}
Let ${\mathcal A}_s$ be a free abelian group with generators
$b(m), m \in 2{\bf N}+1$. 
Define a binary operation $\times: {\mathcal A}_s \times {\mathcal A}_s
\rightarrow {\mathcal A}_s$,
$$b(q) \times b(r)=\sum_{j \in {\bf N}} {\rm dim} \ I \ 
{L(3/2,h_{1,j}) \choose L(3/2,h_{1,q}) \
L(3/2,h_{1,r})} b(j).$$
Then ${\mathcal A}_s$ is a commutative associative
ring, and the mapping $b(m) \mapsto V(\frac{m-1}{4})$
gives an isomorphism to the Grothendieck ring
${\mathcal R}ep({\mathfrak osp}(1|2))$.
\end{corollary}
{\em Proof:} The proof follows from
Theorem \ref{last}(c) and (\ref{osp12}).
\epf

\section{Multiplicity--$2$ fusion rules and super logarithmic intertwiners}

\subsection{Multiplicity--$2$}
We have seen that in the $c=\frac{3}{2}$ case 
all fusion coefficients are $0$ or $1$.
Still, we expect (according to \cite{HM}) that for some
vertex operator superalgebras $L(c,0)$, fusion coefficients are $2$.

Here is one example.
If $c=0$, as in the case of the Virasoro algebra, the super vertex
operator algebra $L(0,0)=\frac{M(0,0)}{\langle G(-1/2)v_0,G(-3/2)v_0
\rangle}$ is trivial. 
Still we can consider a vertex operator superalgebra
$V(0,0):=\frac{M(0,0)}{\langle G(-1/2)v \rangle}$
Clearly, for every $h \in \mathbb{C}$,
we have (all modules are considered to be $V(0,0)$--modules):
\begin{equation}
{\rm dim} \ I \ {L(0,0) \choose L(0,h) \ L(0,h)}=2.
\end{equation}

The previous example is little bit awkward. 
Here is a nice example with ``irrational'' central charge:
\begin{proposition} \label{twons}
\begin{equation} \label{sqrt}
{\rm dim} \ I \
{L(\frac{15}{2}-3\sqrt{5},\frac{\sqrt{5}}{2}-1) \choose
L(\frac{15}{2}-3\sqrt{5},\frac{3}{4}(\frac{\sqrt{5}}{2}-1)) \
L(\frac{15}{2}-3\sqrt{5},\frac{3}{4}(\frac{\sqrt{5}}{2}-1))}=2.
\end{equation}
\end{proposition}
{\em Proof:}
It is not hard to see (by using a result form \cite{Aa} or \cite{D}) that 
$M(\frac{15}{2}-3\sqrt{5},\frac{3}{4}(\frac{\sqrt{5}}{2}-1))$
has the unique submodule that is irreducible (the case $II_+$ in
\cite{Aa}). If we analyze the determinant formula \cite{KWa}, singular vectors, 
and then use Theorem \ref{main}, after some calculation we obtain
(\ref{sqrt}).
\epf

\subsection{Logarithmic intertwiners}
In \cite{M2} we introduced and constructed several examples of
logarithmic intertwining operators. Roughly, logarithmic intertwiners
exist if matrix coefficients yield some logarithmic solutions.

By straightforward super--extension we obtain the following result:
\begin{proposition} \label{logns}
\begin{equation}
{\rm dim} \ I \ {W_2(\frac{27}{2},\frac{-3}{2}) \choose
L(\frac{27}{2},\frac{-3}{2}) \ L(\frac{27}{2},\frac{-3}{2})}=2
\end{equation}
\end{proposition}
{\em Proof:}
Again, the result follows by combining techniques from this paper
and \cite{M2}
\epf

\section{Future work and open problems}
\begin{itemize}
\item
For which triples $L(c,h_1)$, $L(c,h_2)$ and $L(c,h_3)$ do
we have
$${\rm dim} \ I \ {L(c,h_3) \choose L(c,h_1) \ L(c,h_2)}=2 ?$$
\item
Determine the fusion ring for degenerate
minimal models for $N=2$ superconformal algebra (cf. \cite{M3}).
\item
Construct an analogue of the vertex tensor categories 
constructed in \cite{HM} (by using the main result in \cite{Ad}), for the
models studied in this paper.

\end{itemize}

\noindent {\small \sc Department of Mathematics, Rutgers University,
110 Frelinghuysen Rd., Piscataway, NJ 08854-8019}
\vskip 10mm 

\noindent {\em E-mail address}: amilas@math.rutgers.edu

\end{document}